\documentclass[11pt]{article}
\usepackage{epsf}
\usepackage{psfrag}
\usepackage{amsmath}
\usepackage{amssymb}
\usepackage{lscape}
\usepackage{latexsym}
\usepackage[usenames]{color}
\usepackage[mathcal]{eucal}
\usepackage{stmaryrd}
\usepackage{textcomp}
\usepackage{setspace} % The package are used for spacing.
\usepackage{graphicx} % The package are used for 'jpg ..' figure.
\usepackage{pgfplots} % The package are used for 'pdf' figure.
\usepackage{epsfig}   % The package are used for 'eps' figure.
\usepackage{epstopdf} % The package are used for 'eps' figure.
\usepackage{appendix}
\usepackage{afterpage}
\usepackage{float} % for position setting of figures
\usepackage{array} % column formatting in the comparison table
\usepackage[misc,geometry]{ifsym} % corresponding author '\Letter'
\usepackage[font=small,skip=0pt]{caption}
\newtheorem{lemma}{Lemma}[section]

\newtheorem{proposition}{Proposition}[section]

\newtheorem{remark}{Remark}[section]
\newtheorem{assumption}{Assumption}[section]

\newtheorem{algorithm}{Algorithm}[section]
\newcommand{\ds}{\displaystyle}
\newcommand{\dZ}{{\cal Z \kern -0.7em Z}}
\newcommand{\dC}{{\rm\hbox{C \kern-0.8em\raise0.2ex\hbox{\vrule height5.4pt width0.7pt}}}}
\newcommand{\dQ}{{\rm\hbox{Q \kern-0.85em\raise0.25ex\hbox{\vrule height5.4pt width0.7pt}}}}

\newcommand{\proofbox}{\hspace{\fill}{$\Box$}}

\newcommand\old[1]{}
\newcommand{\beqa}{\begin{eqnarray*}}
\newcommand{\eeqa}{\end{eqnarray*}}
\title{\bf \Large{Minimum-Time Interceptor Selection and Obstacle-Avoiding Path Planning for Heterogeneous UAVs}}
\author{
M. M. Rizvi  \thanks{Federation University Australia, {\tt email:  rizmm001@mymail.unisa.edu.au}, } \and Masuda Akter  \Letter \thanks{Department of Mathematics, University of Chittagong, Chittagong-4331, Bangladesh,
   {\tt email: masuda@cu.ac.bd}.(corresponding author)}}
\begin{document}

\maketitle

\begin{abstract}
\noindent Selecting the most suitable unmanned aerial vehicle (UAV) to intercept a moving target is challenging when multiple vehicles have different speeds, initial headings, turning capabilities and launch delays, particularly in environments containing obstacles. This paper develops a computationally efficient framework that integrates minimum-time interceptor selection with obstacle-avoiding path planning for heterogeneous UAVs. The interception problem for each UAV is formulated as a minimum-time optimal control problem and solved numerically using an arc-parameterised nonlinear programming model based on Dubins-type paths. The candidate trajectories consist of minimum-radius circular arcs, straight-line segments and tangential obstacle-boundary arcs, with exact geometric distance constraints imposed to ensure obstacle avoidance along the entire path. To improve computational efficiency, a closed-form lower bound on interception time is incorporated into a bound-and-select algorithm that eliminates UAVs unable to improve the best verified feasible interception time. The proposed pruning strategy preserves the best computed interception decision obtained through full enumeration while reducing computational effort. Numerical experiments investigate the effects of target speed, UAV turning radius, launch delay and obstacle configuration. Across 20 randomly generated scenarios involving up to eight UAVs, the proposed approach reduces average computational time by approximately 85$\%$ for eight-vehicle teams without changing the selected interceptor or its computed interception time. The results demonstrate that the nearest or fastest UAV is not necessarily the most effective interceptor. An illustrative geographical case study further demonstrates the applicability of the framework to UAV interception planning. The proposed approach provides an efficient method for selecting a single UAV and generating feasible Dubins-type interception trajectories in obstacle-constrained environments.
\end{abstract}


\textbf{AMS} subject classifications: 49M37, 90C30, 49K15, 90C27.

%=======================================================================
\section{Introduction}\label{sec:intro}
%=======================================================================
\noindent Unmanned aerial vehicles (UAVs) are increasingly deployed in teams for surveillance, border protection, search and rescue, and the defence of critical infrastructure. In many of these missions a moving target is detected while several vehicles are available to respond, some already airborne on patrol and others waiting at launch sites. These vehicles usually differ in position, heading, speed, minimum turning radius and launch delay, so they can reach the same target at very different times. The operator, or an onboard decision module, must decide which vehicle to send and along which path, and this decision has to account for the turning limits of fixed-wing vehicles, for no-fly zones and other obstacles, for the time needed to launch a vehicle, and for the motion of the target. These factors interact. A vehicle that faces away from the target must first complete a turn, an obstacle in its way forces a detour, and a vehicle on the ground cannot move until it is launched, while the target keeps moving throughout. As a result, the nearest or the fastest vehicle is not necessarily the one that arrives first.

\noindent A reliable choice therefore requires each vehicle’s actual minimum interception time among the obstacles. For one vehicle, this time is the solution of a minimum-time optimal control problem with a bounded turning rate and obstacle-avoidance constraints, which leads to a nonconvex nonlinear programme. Solving this problem for every vehicle of a team is expensive, whereas the decision has to be made shortly after the target is detected. In addition, a launch decision must rest on a path that is guaranteed to avoid every obstacle.

\noindent In this paper we study the case in which a single vehicle is launched, and we ask how to identify, with certainty about the computed quantities and at low computational cost, the vehicle that intercepts the target first. For each vehicle, the interception problem is solved over a family of candidate path types with exact obstacle-avoidance constraints, and only paths verified to be feasible are accepted. A closed-form lower bound on each vehicle’s interception time, obtained by ignoring the turning limit and the obstacles, is then combined with these verified times in a bound-and-select algorithm. A vehicle is discarded only when its lower bound shows that it cannot beat the best verified time already found. The method therefore returns the same vehicle and the same interception time as solving every vehicle, together with lower and upper bounds on the optimal interception time, at a much lower computational cost.
%\color{black}

\subsection{Literature review}\label{sec:lit}
\noindent \textit{Curvature-constrained paths.} The shortest path of bounded curvature between two configurations in the plane was characterised by Dubins \cite{Dubins1957}: it consists of at most three pieces, each a circular arc of minimum turning radius ($C$) or a straight segment ($S$), of type $CSC$ or $CCC$. Kaya \cite{Kaya2017} revisited this Markov--Dubins problem through optimal control theory, obtaining the optimal path types from the necessary conditions of optimality and computing the paths numerically, and later extended the approach to Markov--Dubins interpolating curves, the shortest curvature-constrained curves that pass through a prescribed sequence of points \cite{Kaya2019}.  When the terminal heading is free, as in interception, the optimal path to a point is of type $CS$ or $CC$ \cite{Bui1994}. These results underlie the planning of fixed-wing UAV paths \cite{Tsourdos2011}, and the Dubins travelling salesperson problem, where a curvature-constrained vehicle visits a set of points, has been studied through lower bounds, heuristics and approximation algorithms \cite{Savla2008,Isaiah2015,Manyam2018}. For the related problem of a UAV touring several targets, Forkan, Rizvi and Chowdhury \cite{Forkan2022} computed optimal paths by a geometric approach and by an arc-parameterisation approach.

\noindent \textit{Interception of moving targets.} When the goal point moves, the terminal condition couples the path and the time. An early study of this problem is due to Looker \cite{Looker2008}. Meyer et al.~\cite{Meyer2015} studied Dubins paths that intercept a target moving at constant velocity, Zheng et al.~\cite{Zheng2021} derived time-optimal guidance for Dubins vehicles intercepting moving targets, Buzikov and Galyaev \cite{Buzikov2021} analysed the time-optimal interception of a moving target by a Dubins car, and Manyam et al.~\cite{Manyam2019} considered a target moving on a circle. These works do not include obstacles.

\noindent \textit{Bounded-curvature paths among obstacles.} With obstacles, the geometry of shortest bounded-curvature paths is considerably harder. Boissonnat and Lazard \cite{Boissonnat1996} gave a polynomial-time algorithm for shortest bounded-curvature paths amidst \emph{moderate} obstacles, whose boundaries have curvature not exceeding that of the vehicle; this is exactly the situation of Assumption (A3) below, where every obstacle circle is at least as large as every pursuer's minimum turning circle. Jha, Chen and Shima \cite{Jha2020} characterised shortest Dubins paths constrained to pass through a circular boundary, and the shortest Dubins path avoiding several circular obstacles has recently been studied in \cite{Li2026}. These works concern a single vehicle and a fixed goal.

\noindent In \cite{Akter2023}, the authors optimised Dubins paths that intercept a moving target in free space, and in \cite{Akter2025} the authors formulated the interception of a moving target by one pursuer that tours a finite number of static circular obstacles, through an arc-parameterised nonlinear programme in which the obstacles are entered and left tangentially. 

\noindent \textit{Multi-vehicle assignment and cooperative interception.} When several vehicles are available, path planning is combined with assignment. The weapon--target assignment problem has a long history in operations research; Ahuja et al.~\cite{Ahuja2007} gave exact and heuristic algorithms for it and Kline et al.~\cite{Kline2019} surveyed the field. For curvature-constrained vehicles, Rathinam et al.~\cite{Rathinam2007} gave an approximation algorithm for the allocation of targets to nonholonomic vehicles, Beard et al.~\cite{Beard2002} coordinated target assignment and interception for teams of UAVs, Shanmugavel, Tsourdos and White \cite{Tsourdos2011} planned cooperative Dubins paths with simultaneous arrival, and Wu et al.~\cite{Wu2023} proposed hierarchical task assignment for large-scale group-to-group interception. Multiple-pursuer single-evader differential games have been analysed geometrically by Von Moll et al.~\cite{VonMoll2019}. In these works the assignment cost is usually an approximation, such as a Euclidean distance, an obstacle-free Dubins length or a heuristic travel time, and the obstacles are either ignored or handled by the lower path-planning layer.
%\color{blue}
\subsection{Research gaps}

\noindent Table~\ref{tab:gaps} compares the present work with the main groups of studies reviewed in Section~1.1. Each row is a feature needed to select, from a heterogeneous team, the vehicle that intercepts a moving target first among obstacles; the gap to which it belongs is given in brackets.

\begin{table}[H]
\centering
\caption{Comparison of the present work for showing the research gap.}
\label{tab:gaps}
%\small
\scriptsize
\renewcommand{\arraystretch}{1.25}
\begin{tabular}{|>{\raggedright\arraybackslash}p{\dimexpr0.168\textwidth-2\tabcolsep}|>{\raggedright\arraybackslash}p{\dimexpr0.159\textwidth-2\tabcolsep}|>{\raggedright\arraybackslash}p{\dimexpr0.159\textwidth-2\tabcolsep}|>{\raggedright\arraybackslash}p{\dimexpr0.159\textwidth-2\tabcolsep}|>{\raggedright\arraybackslash}p{\dimexpr0.167\textwidth-2\tabcolsep}|>{\raggedright\arraybackslash}p{\dimexpr0.188\textwidth-2\tabcolsep}|}
\hline
\textbf{Feature (gap)} & \textbf{Free-space interception} \cite{Meyer2015,Zheng2021,Buzikov2021,Manyam2019,Akter2023} & \textbf{Paths among obstacles} \cite{Boissonnat1996,Jha2020,Li2026} & \textbf{Single-UAV obstacle-avoidance model} \cite{Akter2025} & \textbf{Team assignment} \cite{Ahuja2007,Beard2002,Rathinam2007,Wu2023} & \textbf{Present work} \\
\hline
\textbf{Moving target} & Yes; constant velocity (\cite{Manyam2019}: on a circle) & No; fixed goal & Yes; constant velocity & Mostly fixed targets; moving in \cite{Wu2023} & \textbf{Yes; constant velocity} \\
\hline
\textbf{Turning-radius limit} & Yes & Yes & Yes; one radius & Only in \cite{Beard2002,Rathinam2007} & \textbf{Yes; different radius for each pursuer} \\
\hline
\textbf{Static obstacles} \textit{(G1)} & Not considered & Yes; circular or moderate obstacles & Yes; at most two obstacles toured & Threat regions in \cite{Beard2002}; otherwise not considered & \textbf{Yes; any subset and order of obstacles toured} \\
\hline
\textbf{Number of pursuers} \textit{(G1)} & One & One & One & Several & \textbf{Several, with different speeds and turning radii} \\
\hline
\textbf{Launch delay} \textit{(G1)} & Not considered & Not considered & Not considered & Not considered & \textbf{Yes; different for each pursuer} \\
\hline
\textbf{How the vehicle is chosen} \textit{(G1, G2)} & No choice (one pursuer) & No choice (one vehicle) & No choice (one pursuer) & Approximate cost: distance, obstacle-free path length or heuristic travel time & \textbf{Verified minimum interception time; compared with nearest, fastest, shortest-path and straight-line-time rules} \\
\hline
\textbf{Obstacle clearance} \textit{(G3)} & Not applicable & Exact (geometric) & Checked at sampled points after optimisation & Not exact; left to the path-planning layer & \textbf{Exact distances inside the optimisation; every path verified} \\
\hline
\textbf{Lower bound and pruning} \textit{(G4)} & Not considered & Not considered & Not considered & Bounds for the assignment only \cite{Ahuja2007}; none on paths & \textbf{Closed-form lower bound; certified pruning; bracket $\boldsymbol{T_{\mathrm{LB}} \le T^{*} \le T_{\mathrm{feas}}}$} \\
\hline
\end{tabular}
\end{table}

\noindent As Table~\ref{tab:gaps} shows, earlier studies either solve the interception of a moving target by one vehicle or choose among several vehicles with an approximate cost; none chooses a vehicle by its actual interception time among obstacles with launch delays (G1), compares this choice with practical rules (G2), keeps every path exactly clear of the obstacles (G3), or uses a certified lower bound to avoid solving every vehicle (G4).

\color{black}
  
%\subsection{Research gaps}\label{sec:gaps}
%\noindent The literature leaves the following gaps, which motivate the present work.
%\begin{enumerate}
%\item[(G1)] Interception models for curvature-constrained vehicles either ignore obstacles \cite{Meyer2015,Zheng2021,Buzikov2021,Manyam2019} or consider a single pursuer \cite{Akter2025}; assignment models for teams use approximate travel costs that ignore the turning radius, the obstacles or the launch delay \cite{Ahuja2007,Rathinam2007,Wu2023}.
%\item[(G2)] Selection rules used in practice, such as sending the nearest or the fastest vehicle, have not been compared with the obstacle-aware minimum-time choice for heterogeneous teams.
%\item[(G3)] The single-pursuer problem is a nonconvex nonlinear programme. Methods that solve it numerically rarely state which computed quantities are certified (lower bounds, verified feasible paths) and which are only locally optimal, and collision checks are often performed by sampling the path.
%\item[(G4)] Solving a nonconvex path-planning problem for every vehicle of a team is expensive, and cheap valid bounds that allow vehicles to be discarded without changing the decision have not been exploited for this problem.
%\end{enumerate}
%\newpage
\subsection{Contributions and novelty}\label{sec:contrib}
\noindent The paper offers a computationally efficient, obstacle-aware interceptor-selection framework for heterogeneous curvature-constrained UAVs, with provably valid lower-bound pruning and verified feasible paths. Its contributions are as follows.
\begin{enumerate}
\item[(C1)] \textit{Model.} We extend the single-pursuer obstacle-touring model of \cite{Akter2025} to pursuers with launch delays and individual turning radii, we let every pursuer choose the subset, the order and the touring direction of the obstacles it passes, we enlarge the first block of the path so that both turning directions and $CC$ starts are available, and we impose obstacle avoidance on the whole path, including the obstacles that are not toured, by exact point-to-segment and point-to-arc distances inside the optimisation (Section \ref{sec:single}).
\item[(C2)] \textit{Certified selection.} We combine a closed-form lower bound (Lemma \ref{lem:lb}) with verified feasible interception times in a bound-and-select algorithm and prove that pruning never removes a pursuer that could improve the incumbent, so that the algorithm returns the same pursuer and the same verified time as full enumeration, together with the bracket $T_{\rm LB}\leq T^*\leq T_{\rm feas}$ (Proposition \ref{prop:cert}).
\item[(C3)] \textit{Clear separation of claims.} We state explicitly which results are rigorous (the decomposition of the team problem, the lower bound, the certified pruning and the feasibility of every reported path) and which are heuristic (the optimisation within one pursuer, which uses a restricted path family and multistart local optimisation).
\item[(C4)] \textit{Evidence.} We report solver diagnostics for every solve, repeated runs with different random seeds, sensitivity analyses in four parameters, a comparison with three practical selection rules on random scenarios, a computational comparison with full enumeration for team sizes up to eight, and a geographic case study. The MATLAB code that reproduces all tables and figures uses exactly the model and the data of the paper.
\end{enumerate}

\noindent The paper is organised as follows. Section \ref{sec:problem} states the problem and assumptions. Section \ref{sec:single} presents the single-pursuer model and the candidate path family. Section \ref{sec:selection} presents the team problem, the lower bound and the certified bound-and-select algorithm. Section \ref{sec:numerics} reports the numerical results, the sensitivity analyses, the comparisons and the computational study, and Section \ref{sec:case} the case study. Section \ref{sec:discussion} discusses limitations and Section \ref{sec:conclusion} concludes.

%=======================================================================
\section{Problem Statement}\label{sec:problem}
%=======================================================================
\noindent Consider $K$ pursuers $P_1,\dots,P_K$ moving in the plane. Pursuer $P_k$ starts from the configuration $(x_{k,0},y_{k,0},\theta_{k,0})$, moves forward at constant speed $|\mathbf{V}_k|$, and has minimum turning radius $\rho_k$, so that its curvature is bounded by $a_k=1/\rho_k$. It can begin to move only after a launch delay $\delta_k\geq 0$, measured from the detection time $t_0=0$; $\delta_k$ models the time needed to prepare and launch a vehicle on the ground, and $\delta_k=0$ for a vehicle that is already airborne (in which case $(x_{k,0},y_{k,0},\theta_{k,0})$ is its configuration at $t_0$). The target starts from $(x_{T,0},y_{T,0})$ and moves along a straight line with constant speed $|\mathbf{V}_T|$ and heading $\theta_T$,
\begin{equation}\label{eq:target}
x_T(t)=x_{T,0}+|\mathbf{V}_T|\cos\theta_T\, t, \qquad y_T(t)=y_{T,0}+|\mathbf{V}_T|\sin\theta_T\, t .
\end{equation}
The environment contains $n$ static obstacles. Obstacle $O_j$ is enclosed by a circle, called its turning circle, with centre $(x_j,y_j)$ and radius $R_j$, which no pursuer may enter. The task is to choose one pursuer and a path for it so that the target is intercepted as early as possible. We make the following assumptions.
\begin{assumption}\label{ass:main}
\begin{enumerate}
\item[(A1)] Every pursuer is faster than the target, $|\mathbf{V}_k|>|\mathbf{V}_T|$, $k=1,\dots,K$.
\item[(A2)] The target's initial state, speed and heading are known at $t_0$ and do not change during the engagement.
\item[(A3)] Every obstacle turning circle is at least as large as the minimum turning circle of every pursuer, $R_j\geq\rho_k$ for all $j$ and $k$, and the obstacle turning circles are pairwise disjoint.
\item[(A4)] The pursuers act independently: no constraint couples their paths, and interception occurs when a launched pursuer reaches the target.
\item[(A5)] No pursuer starts inside an obstacle turning circle.
\end{enumerate}
\end{assumption}
\noindent Assumption (A3) guarantees that every pursuer can fly along the boundary of every obstacle turning circle with an admissible curvature $1/R_j\leq a_k$; in the terminology of \cite{Boissonnat1996}, the obstacles are moderate for every pursuer. All scenarios in this paper satisfy (A1)--(A5), and the MATLAB code checks (A1), (A3) and (A5) before solving.

%=======================================================================
\section{Minimum-Time Interception by One Pursuer}\label{sec:single}
%=======================================================================
\subsection{Optimal control formulation}
\noindent For pursuer $P_k$, let $(x_k(t),y_k(t),\theta_k(t))$ denote its position and heading. Its interception problem is
\begin{equation}\label{eq:ocp}
\begin{array}{cll}
\ds T_k^*=\min_{u_k(\cdot),\;t_f} & t_f & \\[3mm]
\mbox{subject to} & \dot{x}_k(t)=|\mathbf{V}_k|\cos\theta_k(t),\quad \dot{y}_k(t)=|\mathbf{V}_k|\sin\theta_k(t), & t\in[\delta_k,t_f],\\[2mm]
& \dot{\theta}_k(t)=|\mathbf{V}_k|\,u_k(t),\qquad |u_k(t)|\leq a_k, & t\in[\delta_k,t_f],\\[2mm]
& \big(x_k(\delta_k),y_k(\delta_k),\theta_k(\delta_k)\big)=(x_{k,0},y_{k,0},\theta_{k,0}), & \\[2mm]
& \big(x_k(t)-x_j\big)^2+\big(y_k(t)-y_j\big)^2\geq R_j^2, & t\in[\delta_k,t_f],\ j=1,\dots,n,\\[2mm]
& x_k(t_f)=x_T(t_f),\qquad y_k(t_f)=y_T(t_f), &
\end{array}
\end{equation}
where $x_T,y_T$ are given by \eqref{eq:target}. The decision variables are the curvature control $u_k(\cdot)$ and the interception time $t_f$; the terminal heading is not constrained, because only the positions have to coincide. If \eqref{eq:ocp} is infeasible we set $T_k^*=+\infty$. Since the speed is constant, the path length of $P_k$ is $|\mathbf{V}_k|(t_f-\delta_k)$, so that for one pursuer minimising $t_f$ is equivalent to minimising the path length.

\subsection{Structure of optimal paths and the candidate path family}\label{sec:family}
\noindent By Pontryagin's maximum principle \cite{Pontryagin1962}, on subarcs where no obstacle constraint is active the optimal control satisfies $u_k\in\{-a_k,0,a_k\}$ (bang--bang or singular), so that the path consists of arcs of radius $\rho_k$ and straight segments. On a boundary subarc the pursuer follows an obstacle turning circle with curvature $\pm1/R_j$, which is admissible by (A3), and the path enters and leaves the circle tangentially, since a corner is not admissible for a curvature-bounded path \cite{Boissonnat1996,Jha2020,Akter2025}. A free subarc between two boundary subarcs connects two configurations whose headings are tangent to two disjoint circles of radius at least $\rho_k$; a common tangent segment of the two circles is then the natural candidate, and it is the only connection used in our family.

\noindent We therefore parameterise a path of $P_k$ by a \emph{path type} $\pi=(\phi,\,(j_1,\dots,j_m),\,(\sigma_1,\dots,\sigma_m))$, where
\begin{itemize}
\item $\phi\in\{LR,RL\}$ is the order of the two arcs on the pursuer's minimum turning circle at the start of the path. Both arc lengths are free and may be zero, so the start of the path can be a left turn ($L$), a right turn ($R$), or a turn followed by a turn in the opposite direction ($LR$ or $RL$, i.e. a $CC$ start), followed by a straight segment;
\item $(j_1,\dots,j_m)$ is the ordered list of the $m\geq 0$ obstacles that the pursuer tours, with distinct indices;
\item $\sigma_b\in\{+1,-1\}$ is the touring direction of $O_{j_b}$: $\sigma_b=+1$ if $O_{j_b}$ lies on the pursuer's left and is toured counterclockwise, and $\sigma_b=-1$ if it lies on the right and is toured clockwise.
\end{itemize}
The path is divided into $m+1$ blocks. Block $1$ consists of two arcs on the minimum turning circle, in the order $\phi$, and a straight segment; block $b+1$ ($b=1,\dots,m$) consists of an arc on the turning circle of $O_{j_b}$ in the direction $\sigma_b$ and a straight segment. The straight segment of block $b\leq m$ ends tangentially on the turning circle of $O_{j_b}$, and the straight segment of block $m+1$ ends at the interception point. We write the path type as a word: for example, $RS+LS(O_1)+S_T$ denotes an initial right turn and straight segment, a counterclockwise tour of $O_1$, and a straight segment to the interception point, and $LRS+S_T$ an obstacle-free path that starts with a left and then a right turn. For $n$ obstacles and tours of at most $\bar m$ obstacles, the number of path types is
\begin{equation}\label{eq:count}
N(n,\bar m)=2\sum_{m=0}^{\bar m}\frac{n!}{(n-m)!}\,2^m ,
\end{equation}
so that $N(3,3)=158$ and $N(2,2)=26$. The family $\mathcal{F}$ used in the paper takes $\bar m=n$, i.e. tours of every subset of obstacles in every order are allowed; it contains the family of \cite{Akter2025}, where the start of the path is a single turn and at most two obstacles are toured ($62$ types for $n=3$).

\noindent The family $\mathcal{F}$ is a restriction of the admissible paths of \eqref{eq:ocp}: it does not contain paths that make additional minimum-radius turns between two obstacles or between the last obstacle and the target, or that tour the same obstacle twice. For this reason the value computed over $\mathcal{F}$ is an upper bound on $T_k^*$, and we do not claim that it equals $T_k^*$ in general; Section \ref{sec:family-num} reports how the computed values change when the family is enlarged.

\subsection{Arc-parameterised model}\label{sec:model}
\noindent Let $\ell_1,\dots,\ell_{3m+3}\geq 0$ be the lengths of the pursuer segments and $\ell_T$ the distance travelled by the target until interception. Block $b$ contains the segments $3b-2$, $3b-1$ (arcs) and $3b$ (straight). Its arcs have signed curvatures $(\kappa_{3b-2},\kappa_{3b-1})$, where $(\kappa_1,\kappa_2)=(a_k,-a_k)$ if $\phi=LR$ and $(-a_k,a_k)$ if $\phi=RL$, and $(\kappa_{3b+1},\kappa_{3b+2})=(1/R_{j_b},-1/R_{j_b})$ for $b=1,\dots,m$; straight segments have $\kappa_{3b}=0$. With $\tau_i$ the end of segment $i$, $\tau_0=\delta_k$ and $(x_k(\tau_0),y_k(\tau_0),\theta_k(\tau_0))=(x_{k,0},y_{k,0},\theta_{k,0})$, integration of \eqref{eq:ocp} gives, for every segment $i$,
\begin{equation}\label{eq:seg}
\begin{array}{l}
\theta_k(\tau_i)=\theta_k(\tau_{i-1})+\kappa_i\ell_i,\\[2mm]
\ds x_k(\tau_i)=x_k(\tau_{i-1})+\left\{\begin{array}{ll}\ds\frac{\sin\theta_k(\tau_i)-\sin\theta_k(\tau_{i-1})}{\kappa_i}, & \kappa_i\neq0,\\[3mm] \ell_i\cos\theta_k(\tau_{i-1}), & \kappa_i=0,\end{array}\right.\\[8mm]
\ds y_k(\tau_i)=y_k(\tau_{i-1})+\left\{\begin{array}{ll}\ds\frac{\cos\theta_k(\tau_{i-1})-\cos\theta_k(\tau_i)}{\kappa_i}, & \kappa_i\neq0,\\[3mm] \ell_i\sin\theta_k(\tau_{i-1}), & \kappa_i=0.\end{array}\right.
\end{array}
\end{equation}
For a path type $\pi$, the model of pursuer $P_k$ is
\begin{equation}\label{eq:nlp}
\begin{array}{cll}
\ds T_k(\pi)=\min_{\ell_1,\dots,\ell_{3m+3},\,\ell_T} & \ds \delta_k+\frac{1}{|\mathbf{V}_k|}\sum_{i=1}^{3m+3}\ell_i & \\[4mm]
\mbox{subject to} & x_k(\tau_{3b})=x_{j_b}+\sigma_b R_{j_b}\sin\theta_k(\tau_{3b}), & b=1,\dots,m,\\[2mm]
& y_k(\tau_{3b})=y_{j_b}-\sigma_b R_{j_b}\cos\theta_k(\tau_{3b}), & b=1,\dots,m,\\[2mm]
& x_k(\tau_{3m+3})=x_{T,0}+\ell_T\cos\theta_T, & \\[2mm]
& y_k(\tau_{3m+3})=y_{T,0}+\ell_T\sin\theta_T, & \\[2mm]
& \ds \ell_T=|\mathbf{V}_T|\Big(\delta_k+\frac{1}{|\mathbf{V}_k|}\sum_{i=1}^{3m+3}\ell_i\Big), & \\[3mm]
& d_j(\ell)\geq R_j, & j=1,\dots,n,\\[2mm]
& \ell_{3b+2}=0 \ \mbox{if}\ \sigma_b=+1,\qquad \ell_{3b+1}=0 \ \mbox{if}\ \sigma_b=-1, & b=1,\dots,m,\\[2mm]
& 0\leq \ell_{3b-2},\ell_{3b-1}\leq 2\pi/|\kappa_{3b-2}|,\quad \ell_{3b}\geq 0, & b=1,\dots,m+1,\\[2mm]
& \ell_T\geq 0, &
\end{array}
\end{equation}
together with \eqref{eq:seg}. The first two constraints place the end of the straight segment of block $b$ on the turning circle of $O_{j_b}$ with the radius orthogonal to the heading, so that the pursuer enters the circle tangentially and the following arc, of curvature $\sigma_b/R_{j_b}$, is centred at $(x_{j_b},y_{j_b})$; the exit is tangential because the next segment is straight. The third to fifth constraints are the interception and time-synchronisation conditions, which include the launch delay: the target moves for $\delta_k$ time units before $P_k$ starts.

\noindent The sixth constraint is the exact obstacle-avoidance condition. For each obstacle $O_j$, $d_j(\ell)$ is the minimum distance from the centre $c_j=(x_j,y_j)$ to the \emph{whole} path, computed segment by segment. For a straight segment from $A$ with unit direction $e$ and length $L$, the distance is $\|A+t^*e-c_j\|$ with $t^*=\min\{\max\{(c_j-A)^\top e,0\},L\}$. For an arc with centre $Q$, radius $r$, start angle $\varphi_0$ and signed sweep $\Delta\varphi$, the distance is $\big|\|c_j-Q\|-r\big|$ if the direction of $c_j-Q$ lies within the swept angle, and the smaller of the distances to the two end points of the arc otherwise. Thus a path satisfies $d_j(\ell)\geq R_j$ for all $j$ if and only if no point of it lies inside an obstacle turning circle; no sampling is involved. The constraint is imposed for every obstacle, including those that are not toured, and it is active (with equality) along the toured arcs. The last constraints select the touring arc that matches $\sigma_b$ and exclude arcs of more than one full turn.

\noindent \textit{Verified feasibility.} Problem \eqref{eq:nlp} is a nonconvex nonlinear programme, solved by a local method from several starting points. Every local solution is checked a posteriori: it is \emph{verified feasible} if the largest residual of the equality constraints is below $\varepsilon_{\rm eq}D_k$ and the exact clearance $\min_j\{d_j(\ell)-R_j\}$ is at least $-\varepsilon_{\rm cl}D_k$, where $D_k$ is the initial pursuer--target distance and $\varepsilon_{\rm eq}=10^{-7}$, $\varepsilon_{\rm cl}=10^{-6}$. A verified feasible solution defines a path that is admissible for \eqref{eq:ocp} (up to these tolerances) and intercepts the target at the reported time. We keep \emph{every} verified feasible local solution of every path type: a path type is never discarded because one of its local solutions collides, and the value of a type is the smallest objective over its verified feasible runs. The value of pursuer $P_k$ is
\begin{equation}\label{eq:That}
\widehat T_k=\min_{\pi\in\mathcal{F}}\ \widehat T_k(\pi),\qquad \widehat T_k(\pi)=\min\{\mbox{objective of the verified feasible runs of type }\pi\},
\end{equation}
with $\widehat T_k(\pi)=+\infty$ if no run of type $\pi$ is verified feasible. Since every verified feasible run is an admissible path of \eqref{eq:ocp}, $\widehat T_k\geq T_k^*$.

\noindent The structural observations of Section \ref{sec:selection} (decomposition of the team problem and launch of a single pursuer) are elementary consequences of the assumptions and are stated as observations rather than as main results; the main methodological results are the lower bound, the certified pruning and the experimental evidence.
%=======================================================================
\section{Interceptor Selection}\label{sec:selection}
%=======================================================================
\subsection{The team problem and two structural observations}
\noindent Let $\lambda_k\in\{0,1\}$ indicate whether pursuer $P_k$ is launched. Since one pursuer is launched, the team problem is
\begin{equation}\label{eq:team}
\begin{array}{cl}
\ds T^*=\min_{\lambda,\;u_1(\cdot),\dots,u_K(\cdot),\;t_f} & t_f\\[3mm]
\mbox{subject to} & \ds\sum_{k=1}^{K}\lambda_k=1,\quad \lambda_k\in\{0,1\},\\[3mm]
& \lambda_k=1\ \Rightarrow\ \big(u_k(\cdot),t_f\big)\ \mbox{is feasible for \eqref{eq:ocp}},\quad k=1,\dots,K.
\end{array}
\end{equation}

\begin{remark}[Decomposition]\label{obs:decomp}
Under Assumption \ref{ass:main}, $T^*=\min_{k}T_k^*$, and launching a pursuer $k^*$ with $T_{k^*}^*=T^*$ along an optimal path of \eqref{eq:ocp} solves \eqref{eq:team}. Indeed, for $\lambda=e_k$ the constraints of \eqref{eq:team} reduce to those of \eqref{eq:ocp} for pursuer $k$, because by (A4) no constraint couples the pursuers; the feasible set of \eqref{eq:team} is therefore the disjoint union of the single-pursuer feasible sets, and the minimum over a finite union is the smallest of the minima.
\end{remark}

\begin{remark}[One launch suffices]\label{obs:single}
If a set $S$ of pursuers is launched and interception occurs when the first of them reaches the target, the earliest interception time is $\min_{k\in S}T_k^*\geq T^*$, with equality for $S=\{k^*\}$. This holds because, by (A2) and (A4), the target's trajectory does not depend on which pursuers are launched and the launched pursuers do not interact. Additional launches therefore cannot shorten the interception time; they consume fuel and airframe availability and remove vehicles that could respond to a second threat. If the target's motion were uncertain, a second launch could serve as a backup (Section \ref{sec:discussion}).
\end{remark}

\noindent Remarks \ref{obs:decomp} and \ref{obs:single} are elementary; we state them to make the selection rule precise. Two further consequences matter for the experiments. First, pursuers must be compared by interception time and not by path length, since $T_k^*=\delta_k+L_k^*/|\mathbf{V}_k|$, where $L_k^*$ is the optimal path length, and a pursuer with a shorter path can arrive later if it is slower or has a larger launch delay. Second, neither the distance to the target nor the speed alone identifies the best pursuer: the initial heading, the turning radius, the launch delay and the obstacles between a pursuer and the target all affect the ranking. Section \ref{sec:compare} quantifies how often these simpler rules fail.

\subsection{A closed-form lower bound}
\noindent Solving \eqref{eq:nlp} for all path types of all pursuers is the main computational cost. The following bound is used to skip pursuers that cannot be selected.
\begin{lemma}\label{lem:lb}
Let $\mathbf{d}_k=(x_{T,0}-x_{k,0},\,y_{T,0}-y_{k,0})$, $\mathbf{w}=|\mathbf{V}_T|(\cos\theta_T,\sin\theta_T)$ and
\begin{equation}\label{eq:lbquad}
g_k(t)=\|\mathbf{d}_k+\mathbf{w}t\|^2-|\mathbf{V}_k|^2(t-\delta_k)^2 .
\end{equation}
Under (A1), $g_k$ has two real roots $r_k^-\leq r_k^+$, and $r_k^+\geq\delta_k$. Let $\underline{T}_k=r_k^+$. Then every feasible point $(u_k,t_f)$ of \eqref{eq:ocp} satisfies $t_f\geq\underline{T}_k$; in particular $\underline{T}_k\leq T_k^*$.
\end{lemma}
\noindent \textit{Proof.} The leading coefficient of the quadratic $g_k$ is $|\mathbf{V}_T|^2-|\mathbf{V}_k|^2<0$ by (A1), so $g_k$ is strictly concave and tends to $-\infty$ as $|t|\to\infty$. Since $g_k(\delta_k)=\|\mathbf{d}_k+\mathbf{w}\delta_k\|^2\geq0$, $g_k$ has two real roots $r_k^-\leq\delta_k\leq r_k^+$, and $g_k(t)>0$ for $t\in(r_k^-,r_k^+)$. Let $(u_k,t_f)$ be feasible for \eqref{eq:ocp}. Then $t_f\geq\delta_k$, and the path of $P_k$, of length $|\mathbf{V}_k|(t_f-\delta_k)$, joins the start point to the interception point $(x_T(t_f),y_T(t_f))$, whose distance from the start point is $\|\mathbf{d}_k+\mathbf{w}t_f\|$. Hence $\|\mathbf{d}_k+\mathbf{w}t_f\|\leq|\mathbf{V}_k|(t_f-\delta_k)$, i.e. $g_k(t_f)\leq0$, which together with $t_f\geq\delta_k\geq r_k^-$ forces $t_f\geq r_k^+$ (if $g_k(t_f)=0$ with $t_f=r_k^-$, then $r_k^-=\delta_k$ and $\mathbf{d}_k+\mathbf{w}\delta_k=0$, so $r_k^-=r_k^+=\delta_k$). Taking the infimum over feasible points gives $\underline{T}_k\leq T_k^*$. \proofbox

\noindent The bound $\underline{T}_k$ is the interception time of a vehicle that can turn instantly and ignores the obstacles. It costs one square root, and it is exact when the straight line from the start point to the interception point is admissible, i.e. when the initial heading points at the interception point and the line misses every obstacle. We write $T_{\rm LB}=\min_k\underline{T}_k$.

\subsection{Efficient Interceptor Selection Algorithm}
\noindent Algorithm \ref{alg:select} processes the pursuers in increasing order of their lower bounds and keeps the best verified feasible interception time found so far as the incumbent.

\begin{algorithm}\label{alg:select}\rm (Bound-and-select)
\begin{description}
\item[Step 1] \textbf{(Bounds)} For $k=1,\dots,K$, compute $\underline{T}_k$ from Lemma \ref{lem:lb}, and order the pursuers so that $\underline{T}_{(1)}\leq\dots\leq\underline{T}_{(K)}$. Set the incumbent $\overline{T}=+\infty$.
\item[Step 2] \textbf{(Evaluate or prune)} For $i=1,\dots,K$, with $k=(i)$: if $\underline{T}_k\geq\overline{T}$, prune $P_k$ (and, since the bounds are sorted, all remaining pursuers). Otherwise solve \eqref{eq:nlp} for every path type of $\mathcal{F}$ from several starting points, keep every verified feasible run, and compute $\widehat T_k$ from \eqref{eq:That}. If $\widehat T_k<\overline{T}$, set $\overline{T}=\widehat T_k$ and $k^*=k$.
\item[Step 3] \textbf{(Launch)} Launch $P_{k^*}$ along its verified path; report $T_{\rm feas}=\overline{T}$, the interception point $\big(x_T(T_{\rm feas}),y_T(T_{\rm feas})\big)$ and the bound $[T_{\rm LB},T_{\rm feas}]$.
\end{description}
\end{algorithm}

\noindent The following proposition establishes the validity of the proposed lower-bound pruning rule. It shows that a pursuer can be excluded from further consideration when its lower bound is not smaller than the best verified feasible interception time obtained so far. Consequently, the proposed algorithm achieves the same best computed interception time as full enumeration while reducing computational effort. The proposition also provides lower and upper bounds on the optimal team interception time.

\begin{proposition}\label{prop:cert}
Let Assumption \ref{ass:main} hold and let $\widehat T_k$ be computed by any method that returns verified feasible paths, so that $\widehat T_k\geq T_k^*$ for every $k$. Then:
\begin{enumerate}
\item[(a)] every pursuer pruned in Step 2 satisfies $T_k^*\geq\overline{T}$ at the moment it is pruned; no pruned pursuer can intercept the target earlier than the incumbent path;
\item[(b)] Algorithm \ref{alg:select} returns $T_{\rm feas}=\min_{k}\widehat T_k$, i.e. the same verified time (and, up to ties, the same pursuer) as solving \eqref{eq:nlp} for all pursuers;
\item[(c)] $T_{\rm LB}\leq T^*\leq T_{\rm feas}$, and the returned path is admissible for \eqref{eq:team} and intercepts the target at $T_{\rm feas}$.
\end{enumerate}
\end{proposition}
\noindent \textit{Proof.} (a) If $P_k$ is pruned, then $T_k^*\geq\underline{T}_k\geq\overline{T}$ by Lemma \ref{lem:lb}. (b) The incumbent never increases and every solved pursuer contributes $\widehat T_k$, so $\overline{T}\leq\widehat T_k$ for every solved $k$. For a pruned $k$, $\widehat T_k\geq T_k^*\geq\underline{T}_k\geq\overline{T}$ by (a), so pruned pursuers cannot attain a smaller value than the incumbent, and the final incumbent equals $\min_k\widehat T_k$. (c) The lower inequality follows from Remark \ref{obs:decomp} and Lemma \ref{lem:lb}: $T^*=\min_kT_k^*\geq\min_k\underline{T}_k$. The upper inequality holds because the returned path is verified feasible for \eqref{eq:ocp} with $k=k^*$, hence, by Remark \ref{obs:decomp}, feasible for \eqref{eq:team}. \proofbox

\noindent Proposition \ref{prop:cert} separates the certified and the heuristic parts of the method. The lower bound, the pruning, the feasibility of the returned path and the bracket are rigorous. The optimisation within one pursuer is heuristic in two respects, the restriction to the family $\mathcal{F}$ and the use of multistart local optimisation, so $\widehat T_k$ is an upper bound on $T_k^*$ and $T_{\rm feas}$ need not equal $T^*$; the relative gap $(T_{\rm feas}-T_{\rm LB})/T_{\rm feas}$ is reported for every instance. Pruning is decision-neutral with respect to the computed values: it never changes the selected pursuer or its time, it only avoids work. Note also that the pursuer with the smallest lower bound need not be selected, since $\underline{T}_k$ ignores the turning radius and the obstacles; the bound discards pursuers, it does not select them.

\noindent \textit{Complexity.} Step 1 costs $O(K\log K)$. Each solved pursuer requires $N(n,\bar m)$ multistart solves of \eqref{eq:nlp}, each with $3m+4$ variables, $2m+3$ equality constraints and $n$ exact distance constraints whose evaluation costs $O(mn)$. The effort is dominated by the number of solved pursuers, which Algorithm \ref{alg:select} reduces; Section \ref{sec:comp} measures this reduction.
%=======================================================================
\section{Numerical Results}\label{sec:numerics}
%=======================================================================
\subsection{Implementation and experimental protocol}\label{sec:setup}
\noindent For every path type, model \eqref{eq:nlp} is solved by sequential quadratic programming from random starting points: the arc lengths are drawn uniformly from $[0,\pi/|\kappa_i|]$ and the straight segments and $\ell_T$ from $[0,D_k]$. Each start consists of at most two phases. Phase A solves \eqref{eq:nlp} with the equality constraints and bounds only. If the phase-A solution is not verified feasible because it violates the clearance (its equality residual being small), phase B re-solves the problem from that point with the exact obstacle-avoidance constraints, written as one inequality $\mathrm{dist}(c_j,\mbox{segment }i)\geq R_j$ for every segment $i$ and every obstacle $j$ (the pair formed by a touring arc and its own obstacle is omitted, because that arc lies on the obstacle circle by construction). Writing the constraints per segment, rather than through the minimum over the whole path, keeps them continuous with piecewise continuous derivatives. Every run, after phase A or phase B, is verified with the tests of Section \ref{sec:model}, and every verified feasible run is kept. All path types are solved from $4$ starts (stage 1), and the $5$ best types are solved again from $20$ further starts (stage 2). With $n=3$ obstacles, one pursuer thus requires $158\times4+5\times20=732$ multistart runs.

\noindent The numerical results were obtained using MATLAB's fmincon solver with the sequential quadratic programming (SQP) algorithm, implementing the proposed model, candidate path family, feasibility tolerances, and multistart solution protocol. % and the reproducibility statement at the end of the paper describes how the two were cross-checked.
CPU times are recorded for every solve; we record the exit flag, the largest equality residual, the exact clearance, the number of iterations and of function evaluations, the CPU time and the objective value. Every path reported in the paper was additionally re-integrated in closed form and sampled densely: the interception miss distance is below $2\times10^{-11}$ and the sampled clearance agrees with the exact clearance to $10^{-15}$.

\subsection{Benchmark scenario}\label{sec:bench}
\noindent A target is detected at $(40,45)$, moving south-west ($\theta_T=5\pi/4$) at speed $|\mathbf{V}_T|=3$. Three obstacles lie between the pursuers and the target: $O_1$ with centre $(18,21)$ and radius $3.5$, $O_2$ with centre $(30,22)$ and radius $3$, and $O_3$ with centre $(14,36)$ and radius $3$. The three pursuers are listed in Table \ref{tab:data}. Pursuer $P_1$ is an airborne patrol vehicle; $P_2$ is a slow vehicle close to the target but facing away from it and with a launch delay; $P_3$ is the fastest vehicle, with the largest turning radius, a launch delay and an initial heading away from the target. All quantities are nondimensional. Assumption (A3) holds, since $\max_k\rho_k=2.5\leq\min_jR_j=3$. The data in Table \ref{tab:data} provide the parameter information used.

\begin{table}[H]
\caption{\small{\textit{Benchmark scenario: pursuer data, distance to the target's initial position and lower bound $\underline{T}_k$ of Lemma \ref{lem:lb}.}}}
\label{tab:data}
\footnotesize
\vskip 1em
\centering
\begin{tabular}{|c|c|c|c|c|c|c|}
\hline
Pursuer & $(x_{k,0},y_{k,0},\theta_{k,0})$ & $|\mathbf{V}_k|$ & $\rho_k$ & $\delta_k$ & Distance to target & $\underline{T}_k$ \\
\hline
$P_1$ & $(5,5,\pi/4)$ & 6 & 1 & 0 & 53.15 & 5.9122 \\
$P_2$ & $(46,36,-\pi/2)$ & 4 & 2 & 2 & 10.82 & 8.2431 \\
$P_3$ & $(-2,46,\pi)$ & 9 & 2.5 & 2 & 42.01 & 5.6355 \\
\hline
\end{tabular}
\end{table}

\begin{figure}[H]
\centering
\includegraphics[width=0.82\textwidth]{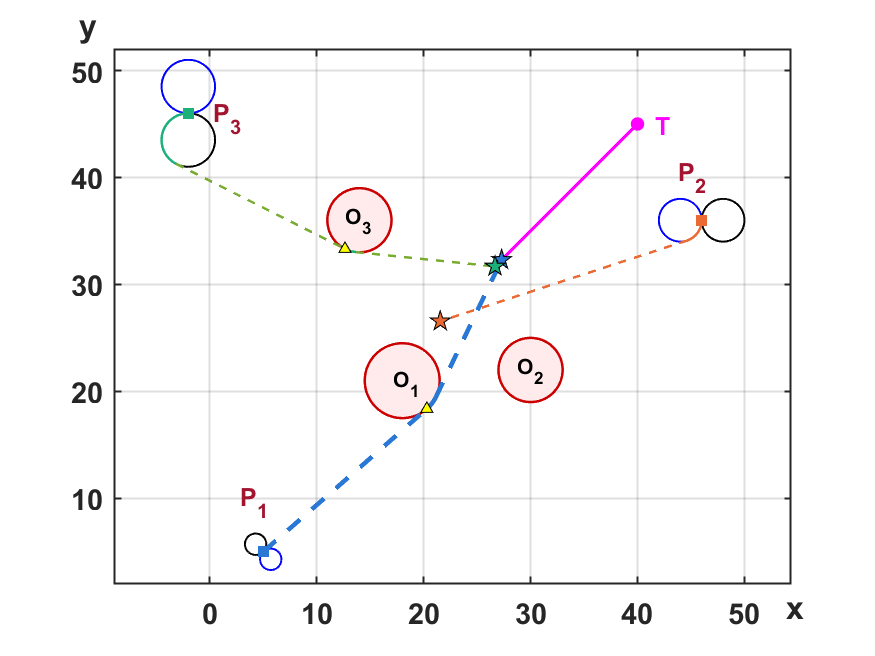}
\caption{\sf\footnotesize Benchmark scenario: best verified interception paths. Only $P_1$ (thick) is launched; the paths of $P_2$ and $P_3$ are those they would follow if launched. Arcs are solid and straight segments dashed; the solid circles are the minimum turning circles at the launch points; triangles mark the tangential entry points; stars mark the interception points.}
\label{fig:bench}
\end{figure}

\noindent Table \ref{tab:results} lists the best verified path of each pursuer and Figure \ref{fig:bench} shows the three paths. Pursuer $P_1$ intercepts the target at $\widehat T_1=5.9898$: after a very short right turn ($\ell_2=0.0676$) it flies straight to the turning circle of $O_1$, tours it counterclockwise for a length $\ell_4=1.4214$, and flies straight to the interception point $(27.29,32.29)$. The straight path from $P_1$ to the target crosses $O_1$, which explains the difference from its lower bound $\underline{T}_1=5.9122$. Pursuer $P_3$ needs $\widehat T_3=6.2779$: it starts with a left turn of length $6.6874$ from its westward heading and then tours $O_3$ counterclockwise. Pursuer $P_2$ needs $\widehat T_2=8.6937$: its obstacle-free path is feasible, but its launch delay and low speed let the target travel $\ell_T=26.08$ before interception. Pursuer $P_1$ is selected, and the team interception time satisfies
\[
T_{\rm LB}=\underline{T}_3=5.6355\ \leq\ T^*\ \leq\ T_{\rm feas}=5.9898 ,
\]
a relative gap of $5.9\%$. For the selected pursuer alone the gap between $\underline{T}_1$ and $\widehat T_1$ is $1.3\%$.

\begin{table}[H]
\caption{\small{\textit{Benchmark scenario: best verified path of each pursuer. $L_k$ is the path length of the pursuer, $\ell_T$ the distance travelled by the target and $\widehat T_k$ the verified interception time; ``clearance'' is the exact clearance $\min_j\{d_j(\ell)-R_j\}$ (zero when an obstacle is toured).}}}
\label{tab:results}
\footnotesize
\vskip 1em
\centering
\resizebox{\textwidth}{!}{%
\begin{tabular}{|c|c|c|c|c|c|c|c|}
\hline
Pursuer & Path type & Nonzero segment lengths & $L_k$ & $\ell_T$ & Interception point & Clearance & $\widehat T_k$ \\
\hline
$P_1$ & $RS+LS(O_1)+S_T$ & $\ell_2=0.0676,\ \ell_3=20.2483,\ \ell_4=1.4214,\ \ell_6=14.2012$ & 35.9385 & 17.9693 & $(27.29,32.29)$ & $0$ & $\mathbf{5.9898}$ \\
$P_2$ & $RS+S_T$ & $\ell_2=2.5096,\ \ell_3=24.2653$ & 26.7749 & 26.0811 & $(21.56,26.56)$ & $3.0991$ & 8.6937 \\
$P_3$ & $LS+LS(O_3)+S_T$ & $\ell_1=6.6874,\ \ell_3=17.6635,\ \ell_4=1.0931,\ \ell_6=13.0571$ & 38.5011 & 18.8337 & $(26.68,31.68)$ & $0$ & 6.2779 \\
\hline
\end{tabular}}
\end{table}

\noindent \textit{Solver diagnostics.} Table \ref{tab:diag} summarises the $732$ runs of each pursuer. Most runs end infeasible because most path types are infeasible or poor for a given pursuer (for example, obstacle-free paths of $P_1$ and $P_3$ cross an obstacle); these runs are recorded but never used. Phase B was invoked in $14$--$33$ runs per pursuer and produced one additional verified feasible run; in this scenario every best path was already found in phase A, with the exact constraints satisfied. The best paths satisfy the equality constraints to $2\times10^{-11}$ and have nonnegative exact clearance.

\begin{table}[H]
\caption{\small{\textit{Benchmark scenario: diagnostics over all runs of each pursuer (stage 1 and stage 2). ``Converged'' counts runs with a successful exit flag of the local solver; ``feasible types'' counts path types with at least one verified feasible run.}}}
\label{tab:diag}
\footnotesize
\vskip 1em
\centering
\resizebox{\textwidth}{!}{%
\begin{tabular}{|c|c|c|c|c|c|c|c|c|}
\hline
Pursuer & Runs & Converged & Verified feasible & Phase B (feasible) & Feasible types & Median iter. & Median $f$-evals & CPU [s] \\
\hline
$P_1$ & 732 & 154 & 140 & 14 (0) & 39 / 158 & 32 & 310 & 77.7 \\
$P_2$ & 732 & 143 & 124 & 22 (1) & 47 / 158 & 28.5 & 297 & 86.4 \\
$P_3$ & 732 & 94 & 61 & 33 (0) & 22 / 158 & 32 & 319 & 77.9 \\
\hline
\end{tabular}}
\end{table}

\noindent \textit{Bound-and-select.} Algorithm \ref{alg:select} orders the pursuers by their lower bounds: $P_3$ ($5.6355$), $P_1$ ($5.9122$), $P_2$ ($8.2431$). Solving $P_3$ gives the incumbent $\overline{T}=6.2779$; solving $P_1$ improves it to $5.9898$; since $\underline{T}_2=8.2431\geq5.9898$, $P_2$ is pruned and its $732$ runs are avoided, a reduction of one third of the CPU time ($155.6$\,s instead of $242.0$\,s). By Proposition \ref{prop:cert}, $P_2$ could not have improved the incumbent. The pursuer with the smallest bound, $P_3$, is not the one selected: its bound ignores its turning radius, its initial heading and the obstacle $O_3$ in its way.

\noindent \textit{Selection criteria.} The selected pursuer $P_1$ is the farthest from the target ($53.15$, against $10.82$ for $P_2$ and $42.01$ for $P_3$) and is not the fastest. The nearest pursuer $P_2$ is the slowest to intercept, the fastest pursuer $P_3$ comes second, and the shortest path ($L_2=26.77$) belongs to $P_2$, so selection by path length would also choose the wrong pursuer. Launching $P_2$ or $P_3$ in addition to $P_1$ would not shorten the interception (Remark \ref{obs:single}).

\subsection{Effect of the candidate path family}\label{sec:family-num}
\noindent To assess the restriction to a finite family, the best value of each pursuer was computed over nested sub-families of $\mathcal{F}$: obstacle-free paths only ($m=0$), tours of at most one obstacle, tours of at most two obstacles, the family of \cite{Akter2025} (single initial turn and $m\leq2$, $62$ types), and the full family ($158$ types, with $CC$ starts and tours of all three obstacles). For all three pursuers the value is the same for every family that contains at least one-obstacle tours ($5.9898$, $8.6937$ and $6.2779$); the obstacle-free family is infeasible for $P_1$ and $P_3$. Among the best path types, those with a $CC$ start (e.g. $RS+LS(O_1)+S_T$ for $P_1$) converge to the same paths with one arc of zero length, and tours of two or three obstacles are longer. Thus, in this scenario, enlarging the family beyond one-obstacle tours did not improve any pursuer; it is nevertheless retained in all experiments, since it costs only additional multistart runs and it guards against configurations in which more complex paths are better.

\subsection{Robustness to the random starting points}\label{sec:seeds}
\noindent The benchmark was solved four more times for every pursuer with different random seeds ($2,\dots,5$), with all other settings unchanged. Table \ref{tab:seeds} shows the results. The verified times agree to within $3\times10^{-11}$ across the five seeds, the same path is found in every case (for $P_1$ and $P_2$ once in an equivalent form with a zero-length arc, e.g. $RS+LS(O_1)+S_T$ with $\ell_1=0$), and the selected pursuer is always $P_1$. Since the gap between the first and the second pursuer ($6.2779-5.9898=0.288$) is many orders of magnitude larger than the variation over seeds, the ranking is not affected by the randomness of the multistart procedure in this scenario.

\begin{table}[H]
\caption{\small{\textit{Benchmark scenario: verified times $\widehat T_k$ over five random seeds.}}}
\label{tab:seeds}
\footnotesize
\vskip 1em
\centering
\resizebox{\textwidth}{!}{%
\begin{tabular}{|c|c|c|c|c|c|}
\hline
Pursuer & Best path type & Minimum over seeds & Maximum over seeds & Range & Selected in \\
\hline
$P_1$ & $RS+LS(O_1)+S_T$ & 5.989757 & 5.989757 & $3.8\times10^{-13}$ & 5 / 5 seeds \\
$P_2$ & $RS+S_T$ & 8.693715 & 8.693715 & $2.4\times10^{-11}$ & 0 / 5 \\
$P_3$ & $LS+LS(O_3)+S_T$ & 6.277898 & 6.277898 & $1.4\times10^{-12}$ & 0 / 5 \\
\hline
\end{tabular}}
\end{table}

\subsection{Sensitivity analysis}\label{sec:sens}
\begin{figure}[H]
\centering
\includegraphics[width=\textwidth]{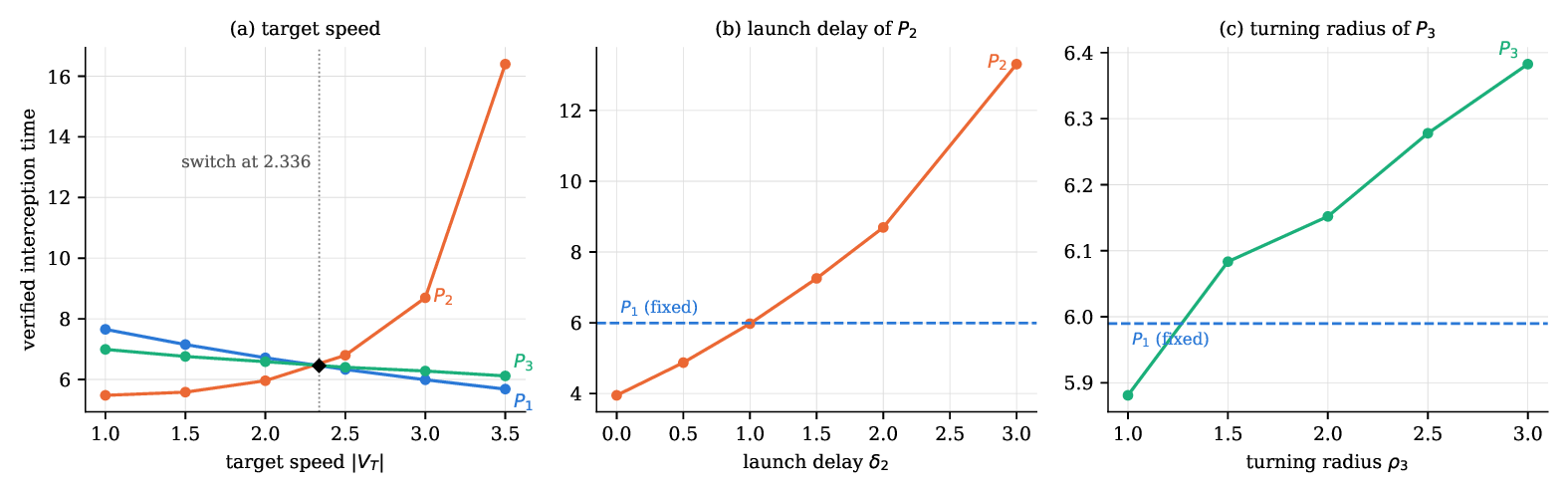}
\caption{\sf\footnotesize Sensitivity of the verified interception times in the benchmark scenario: (a) target speed, with the selection switch from $P_2$ to $P_1$ located by bisection (diamond); (b) launch delay of $P_2$; (c) turning radius of $P_3$. In (b) and (c) the dashed line is the time of $P_1$, which does not depend on the varied parameter.}
\label{fig:sens}
\end{figure}

\noindent We vary one parameter of the benchmark at a time and recompute the affected pursuers with the full protocol of Section \ref{sec:setup}. The results are summarised in Figure \ref{fig:sens} and Tables \ref{tab:sensvt}--\ref{tab:sensobs}.

\noindent \textit{Target speed.} Table \ref{tab:sensvt} lists the verified times of the three pursuers for $|\mathbf{V}_T|\in\{1,1.5,\dots,3.5\}$ (assumption (A1) requires $|\mathbf{V}_T|<4$). A slow target is intercepted first by $P_2$, which starts close to it; as the target speeds up, it moves away from $P_2$ quickly and towards $P_1$ and $P_3$, and from $|\mathbf{V}_T|=2.5$ onwards $P_1$ is selected. The functions $\widehat T_1(|\mathbf{V}_T|)$ and $\widehat T_2(|\mathbf{V}_T|)$ cross between $2.0$ and $2.5$; bisection on $f(v)=\widehat T_2(v)-\widehat T_1(v)$, in which each evaluation re-solves both pursuers with the full protocol, locates the crossover at $|\mathbf{V}_T|^\star\approx2.336$ (final bracket $[2.3359,\,2.3369]$ after nine steps), where both pursuers need $\widehat T\approx6.450$. The best path of $P_3$ changes with the target speed: for slow targets $P_3$ flies an obstacle-free $LS$ or $LRS$ path, while for fast targets it must tour $O_3$. At $|\mathbf{V}_T|=3.5$ the obstacle-free path of $P_2$ is no longer feasible and $P_2$ must tour $O_2$. The lower bounds follow the same trend, but they would have selected $P_3$ for $|\mathbf{V}_T|\geq2.5$, which is never optimal.

\begin{table}[H]
\caption{\small{\textit{Sensitivity to the target speed: verified times (and lower bounds) of the three pursuers; the selected pursuer is in bold.}}}
\label{tab:sensvt}
\footnotesize
\vskip 1em
\centering
\begin{tabular}{|c|c|c|c|c|}
\hline
$|\mathbf{V}_T|$ & $\widehat T_1$ ($\underline{T}_1$) & $\widehat T_2$ ($\underline{T}_2$) & $\widehat T_3$ ($\underline{T}_3$) & Selected \\
\hline
1.0 & 7.6521 (7.5958) & \textbf{5.4762} (4.7307) & 6.9914 (6.2208) & $P_2$ \\
1.5 & 7.1511 (7.0907) & \textbf{5.5809} (4.9676) & 6.7588 (6.0396) & $P_2$ \\
2.0 & 6.7139 (6.6487) & \textbf{5.9588} (5.4384) & 6.5582 (5.8829) & $P_2$ \\
2.5 & \textbf{6.3296} (6.2588) & 6.7976 (6.3344) & 6.4017 (5.7486) & $P_1$ \\
3.0 & \textbf{5.9898} (5.9122) & 8.6937 (8.2431) & 6.2779 (5.6355) & $P_1$ \\
3.5 & \textbf{5.6819} (5.6020) & 16.3962 (14.1086) & 6.1163 (5.5425) & $P_1$ \\
\hline
\end{tabular}
\end{table}

\noindent \textit{Launch delay of $P_2$.} With all other data fixed, $\widehat T_2$ grows from $3.9493$ at $\delta_2=0$ to $13.3092$ at $\delta_2=3$ (Table \ref{tab:sensother}); for $\delta_2=3$ the target has passed $P_2$ and the obstacle-free path is no longer feasible, so $P_2$ must tour $O_1$. Pursuer $P_2$ is selected for $\delta_2\leq1$ ($\widehat T_2=5.9750<5.9898$ at $\delta_2=1$) and $P_1$ for $\delta_2\geq1.5$. A vehicle close to the target is therefore the right choice only if it can be launched quickly; in this scenario one time unit of launch delay costs $P_2$ more than two time units of interception time.

\noindent \textit{Turning radius of $P_3$.} The turning radius of $P_3$ was varied in $[1,3]$, the range allowed by (A3). A smaller radius shortens the initial turn from its westward heading, and $\widehat T_3$ decreases from $6.3826$ ($\rho_3=3$) to $5.8810$ ($\rho_3=1$). For $\rho_3=1$, $P_3$ becomes faster than $P_1$ and is selected; for $\rho_3\geq1.5$, $P_1$ is selected. The lower bound $\underline{T}_3=5.6355$ does not depend on $\rho_3$, which illustrates why it cannot be used to select a pursuer.

\begin{table}[H]
\caption{\small{\textit{Sensitivity to the launch delay of $P_2$ and to the turning radius of $P_3$ (other pursuers unchanged; $\widehat T_1=5.9898$).}}}
\label{tab:sensother}
\footnotesize
\vskip 1em
\centering
\begin{tabular}{|c|c|c|c||c|c|c|c|}
\hline
$\delta_2$ & $\widehat T_2$ & Path type of $P_2$ & Selected & $\rho_3$ & $\widehat T_3$ & Path type of $P_3$ & Selected \\
\hline
0   & 3.9493 & $RS+S_T$ & $P_2$ & 1.0 & 5.8810 & $LS+RS(O_3)+S_T$ & $P_3$ \\
0.5 & 4.8741 & $RS+S_T$ & $P_2$ & 1.5 & 6.0131 & $LS+LS(O_3)+S_T$ & $P_1$ \\
1.0 & 5.9750 & $RS+S_T$ & $P_2$ & 2.0 & 6.1520 & $LS+RS(O_3)+S_T$ & $P_1$ \\
1.5 & 7.2543 & $RS+S_T$ & $P_1$ & 2.5 & 6.2779 & $LS+LS(O_3)+S_T$ & $P_1$ \\
2.0 & 8.6937 & $RS+S_T$ & $P_1$ & 3.0 & 6.3826 & $LS+LS(O_3)+S_T$ & $P_1$ \\
3.0 & 13.3092 & $RS+LS(O_1)+S_T$ & $P_1$ & & & & \\
\hline
\end{tabular}
\end{table}

\noindent \textit{Obstacle configuration.} Table \ref{tab:sensobs} shows the effect of removing each obstacle in turn, removing all of them, and enlarging all radii by $20\%$. Without $O_1$, the straight path of $P_1$ becomes feasible and $\widehat T_1$ drops to its lower bound $5.9122$; without $O_3$, $P_3$ needs no tour ($\widehat T_3=6.2296$); $O_2$ affects no pursuer. With enlarged obstacles all touring paths become longer ($\widehat T_1=6.0255$, $\widehat T_3=6.3000$), while $P_2$ is unaffected. In every configuration $P_1$ remains selected, but the margin over $P_3$ varies between $0.24$ and $0.37$, and in the obstacle-free configuration the optimal path of $P_1$ is a straight line. The obstacles thus change the paths and the interception times; whether they change the selection depends on the margin between the best pursuers, which is examined over random scenarios in Section \ref{sec:compare}.

\begin{table}[H]
\caption{\small{\textit{Sensitivity to the obstacle configuration: verified times and best path types (obstacle labels refer to the original numbering).}}}
\label{tab:sensobs}
\footnotesize
\vskip 1em
\centering
\resizebox{\textwidth}{!}{%
\begin{tabular}{|l|c|c|c|c|}
\hline
Configuration & $\widehat T_1$ (path) & $\widehat T_2$ (path) & $\widehat T_3$ (path) & Selected \\
\hline
All obstacles & \textbf{5.9898} ($RS+LS(O_1)+S_T$) & 8.6937 ($RS+S_T$) & 6.2779 ($LS+LS(O_3)+S_T$) & $P_1$ \\
Without $O_1$ & \textbf{5.9122} ($LS+S_T$) & 8.6937 ($RS+S_T$) & 6.2779 ($LS+LS(O_3)+S_T$) & $P_1$ \\
Without $O_2$ & \textbf{5.9898} ($RS+LS(O_1)+S_T$) & 8.6937 ($RS+S_T$) & 6.2779 ($LS+LS(O_3)+S_T$) & $P_1$ \\
Without $O_3$ & \textbf{5.9898} ($RS+LS(O_1)+S_T$) & 8.6937 ($RS+S_T$) & 6.2296 ($LS+S_T$) & $P_1$ \\
No obstacles & \textbf{5.9122} ($LS+S_T$) & 8.6937 ($RS+S_T$) & 6.2296 ($LS+S_T$) & $P_1$ \\
Radii $\times1.2$ & \textbf{6.0255} ($RS+LS(O_1)+S_T$) & 8.6937 ($RS+S_T$) & 6.3000 ($LS+LS(O_3)+S_T$) & $P_1$ \\
\hline
\end{tabular}}
\end{table}

\subsection{Comparison with practical selection rules on random scenarios}\label{sec:compare}
\noindent To compare the obstacle-aware choice with simpler rules beyond one scenario, we generated $20$ random scenarios in a $60\times60$ region with a fixed seed. Each contains a target with random initial position (in $[25,45]\times[30,50]$), heading and speed ($|\mathbf{V}_T|\in[2,3]$), two non-overlapping obstacles with radii in $[2.5,4]$, and $K=8$ pursuers with random positions and headings, speeds in $[4,9]$, launch delays in $\{0,1,2,3\}$ and turning radii between $1$ and the smallest obstacle radius, so that (A1), (A3) and (A5) hold. %The data are provided in the file \texttt{random\_scenarios.csv}.
Every pursuer of every scenario was solved with the protocol of Section \ref{sec:setup} ($26$ path types, $204$ runs per pursuer, $32{,}640$ runs in total, of which $11{,}131$ were verified feasible), giving the verified time $T_{\rm feas}=\min_k\widehat T_k$ of each scenario. Each selection rule chooses one pursuer, and its \emph{regret} is the relative increase $(\widehat T_{\rm rule}-T_{\rm feas})/T_{\rm feas}$ of the verified interception time. We compare five rules: the nearest pursuer, the fastest pursuer, the pursuer with the shortest verified path, the pursuer with the smallest lower bound $\underline{T}_k$ (straight-line time), and the pursuer with the smallest obstacle-free interception time (model \eqref{eq:nlp} solved without obstacles).

\begin{table}[H]
\caption{\small{\textit{Selection rules on 20 random scenarios with $K=8$ pursuers and $n=2$ obstacles: number of scenarios in which the rule selects a different pursuer from the obstacle-aware choice, and regret in the verified interception time.}}}
\label{tab:rules}
\footnotesize
\vskip 1em
\centering
\begin{tabular}{|l|c|c|c|}
\hline
Selection rule & Wrong choice & Mean regret & Maximum regret \\
\hline
Nearest pursuer & 9 / 20 & 73.5\% & 367.5\% \\
Fastest pursuer & 14 / 20 & 113.3\% & 378.5\% \\
Shortest path length & 9 / 20 & 19.5\% & 115.3\% \\
Smallest lower bound $\underline{T}_k$ & 1 / 20 & 0.3\% & 5.3\% \\
Smallest obstacle-free interception time & 0 / 20 & 0\% & 0\% \\
\hline
\end{tabular}
\end{table}

\noindent Table \ref{tab:rules} shows that the two rules most often used in practice perform poorly: the nearest pursuer is the wrong choice in $9$ of the $20$ scenarios and the fastest in $14$, with mean regrets of $74\%$ and $113\%$ and individual regrets above $360\%$. Selecting by path length, which ignores the differences in speed and launch delay, is wrong in $9$ scenarios. The rules that compare times perform much better: the straight-line time is wrong once (regret $5.3\%$), and the obstacle-free interception time never. In these random scenarios the two obstacles changed the verified time of $24$ of the $160$ pursuers, but never that of the selected pursuer. Random obstacles of this size rarely lie on the best pursuer's route, whereas in the benchmark (Section \ref{sec:bench}) the obstacles lie between all pursuers and the target, and there the straight-line time selects the wrong pursuer ($P_3$). Two conclusions follow. First, the obstacle-aware computation is needed to certify the choice and to provide a feasible path, and it matters most in cluttered environments, where cheap time-based rules can fail. Second, time-based quantities are the right basis for selection, which supports both the time criterion and the use of the lower bound to order the pursuers in Algorithm \ref{alg:select}. The relative gap between $T_{\rm LB}$ and $T_{\rm feas}$ was $9.3\%$ on average (median $5.6\%$, maximum $31.8\%$).

\subsection{Computational benefit of certified pruning}\label{sec:comp}
\begin{figure}[H]
\centering
\includegraphics[width=0.6\textwidth]{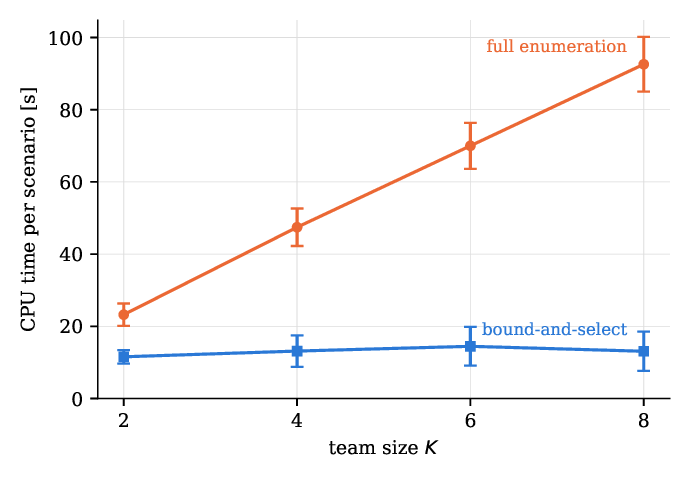}
\caption{\sf\footnotesize CPU time per scenario of full enumeration and of bound-and-select versus the team size (mean $\pm$ one standard deviation over 20 random scenarios).}
\label{fig:cpu}
\end{figure}

\noindent Using the same scenarios, we compared Algorithm \ref{alg:select} with full enumeration for team sizes $K=2,4,6,8$ (the first $K$ pursuers of each scenario). Since the per-pursuer computations are identical in both methods, Algorithm \ref{alg:select} was replayed on the stored per-pursuer results; by Proposition \ref{prop:cert} it returned the same verified time as full enumeration in all $80$ instances, which we checked. Table \ref{tab:cpu} and Figure \ref{fig:cpu} report the fraction of pursuers solved and the CPU time.

\begin{table}[H]
\caption{\small{\textit{Bound-and-select versus full enumeration on the random scenarios (means over 20 scenarios; standard deviations in parentheses). CPU times are serial sums over the solved pursuers.}}}
\label{tab:cpu}
\footnotesize
\vskip 1em
\centering
\begin{tabular}{|c|c|c|c|c|c|}
\hline
$K$ & Pursuers solved & NLP runs solved & CPU full [s] & CPU bound-and-select [s] & Ratio \\
\hline
2 & 50.0\% & 50.0\% & 26.2 (3.2) & 12.7 (2.1) & 0.49 \\
4 & 28.8\% & 28.8\% & 51.4 (4.4) & 14.6 (4.7) & 0.28 \\
6 & 20.8\% & 20.8\% & 76.6 (5.9) & 15.9 (5.9) & 0.21 \\
8 & 14.4\% & 14.4\% & 100.9 (7.1) & 14.6 (5.3) & 0.15 \\
\hline
\end{tabular}
\end{table}

\noindent The effort of full enumeration grows linearly with $K$, about $12.6$\,s per pursuer, while that of bound-and-select remains almost constant: on average only $1.0$, $1.15$, $1.25$ and $1.15$ pursuers had to be solved for $K=2,4,6,8$, so that for $K=8$ the CPU time is reduced by $85\%$ (a factor of seven). In other words, in these scenarios the lower bound identified, after one or two solves, a verified incumbent that no other pursuer could beat. The benefit is largest when the best pursuer is much better than the others, and smallest when several pursuers have similar times. Together with the benchmark (one of three pursuers pruned) and the case study (three of four pruned), this shows that certified pruning removes most of the computational effort without changing the decision.
%=======================================================================
\section{Case Study}\label{sec:case}
%=======================================================================
\noindent We consider an intruding aircraft detected over Feni, Bangladesh, flying towards Dhaka at $650$ km/h. Four UAVs can respond: an airborne patrol vehicle near Satkhira and three vehicles on the ground near Jashore, Dhaka and Chattogram. Three restricted airspaces, centred over Narail, Narayanganj and Cumilla, must not be entered. The geographic coordinates are converted to planar coordinates (in km) with a local equirectangular projection about $(23^\circ\mathrm{N},90^\circ\mathrm{E})$,
\[
x=111.195\,(\lambda-90^\circ)\cos 23^\circ,\qquad y=111.195\,(\varphi-23^\circ),
\]
where $\varphi$ and $\lambda$ are the latitude and longitude in degrees; over the region considered the distance error of this projection is below one percent. Speeds are in km/h, turning radii in km and times in hours. The data are listed in Table \ref{tab:casedata}. The target heads from Feni towards Dhaka, $\theta_T=138.6^\circ$. The vehicle at Dhaka faces west (its runway heading) and needs $3$ min to launch; the vehicle at Jashore is the fastest ($900$ km/h) but needs $6$ min; the vehicle at Chattogram is the nearest to the target ($85$ km) but needs $9$ min and the target flies away from it. Assumption (A3) holds, since $\max_k\rho_k=15\leq\min_jR_j=15$, and all pursuers are faster than the target. The coordinates were taken from public map data and should be read as an illustrative scenario.

\begin{table}[H]
\caption{\small{\textit{Case study: locations, planar coordinates and vehicle data. Obstacles are restricted airspaces of radius $R_j$.}}}
\label{tab:casedata}
\footnotesize
\vskip 1em
\centering
\resizebox{\textwidth}{!}{%
\begin{tabular}{|l|c|c|c|c|c|c|}
\hline
Element & Latitude, longitude & $(x,y)$ [km] & Heading & Speed [km/h] & $\rho_k$ or $R_j$ [km] & Delay [min] \\
\hline
Target (Feni) & $23.0263^\circ$N, $91.3788^\circ$E & $(141.13,\ 2.92)$ & $138.6^\circ$ & 650 & -- & -- \\
$P_1$ Satkhira (airborne) & $22.7377^\circ$N, $89.0571^\circ$E & $(-96.51,\ -29.17)$ & $0^\circ$ & 800 & 10 & 0 \\
$P_2$ Jashore & $23.1664^\circ$N, $89.2081^\circ$E & $(-81.06,\ 18.50)$ & $90^\circ$ & 900 & 15 & 6 \\
$P_3$ Dhaka & $23.8103^\circ$N, $90.4125^\circ$E & $(42.22,\ 90.10)$ & $180^\circ$ & 700 & 8 & 3 \\
$P_4$ Chattogram & $22.3569^\circ$N, $91.7832^\circ$E & $(182.52,\ -71.51)$ & $-90^\circ$ & 750 & 10 & 9 \\
$O_1$ Narail & $23.1956^\circ$N, $89.4971^\circ$E & $(-51.47,\ 21.75)$ & -- & -- & 20 & -- \\
$O_2$ Narayanganj & $23.6432^\circ$N, $90.4822^\circ$E & $(49.36,\ 71.52)$ & -- & -- & 15 & -- \\
$O_3$ Cumilla & $23.4607^\circ$N, $91.1809^\circ$E & $(120.87,\ 51.23)$ & -- & -- & 20 & -- \\
\hline
\end{tabular}}
\end{table}

\begin{figure}[H]
\centering
\includegraphics[width=0.80\textwidth]{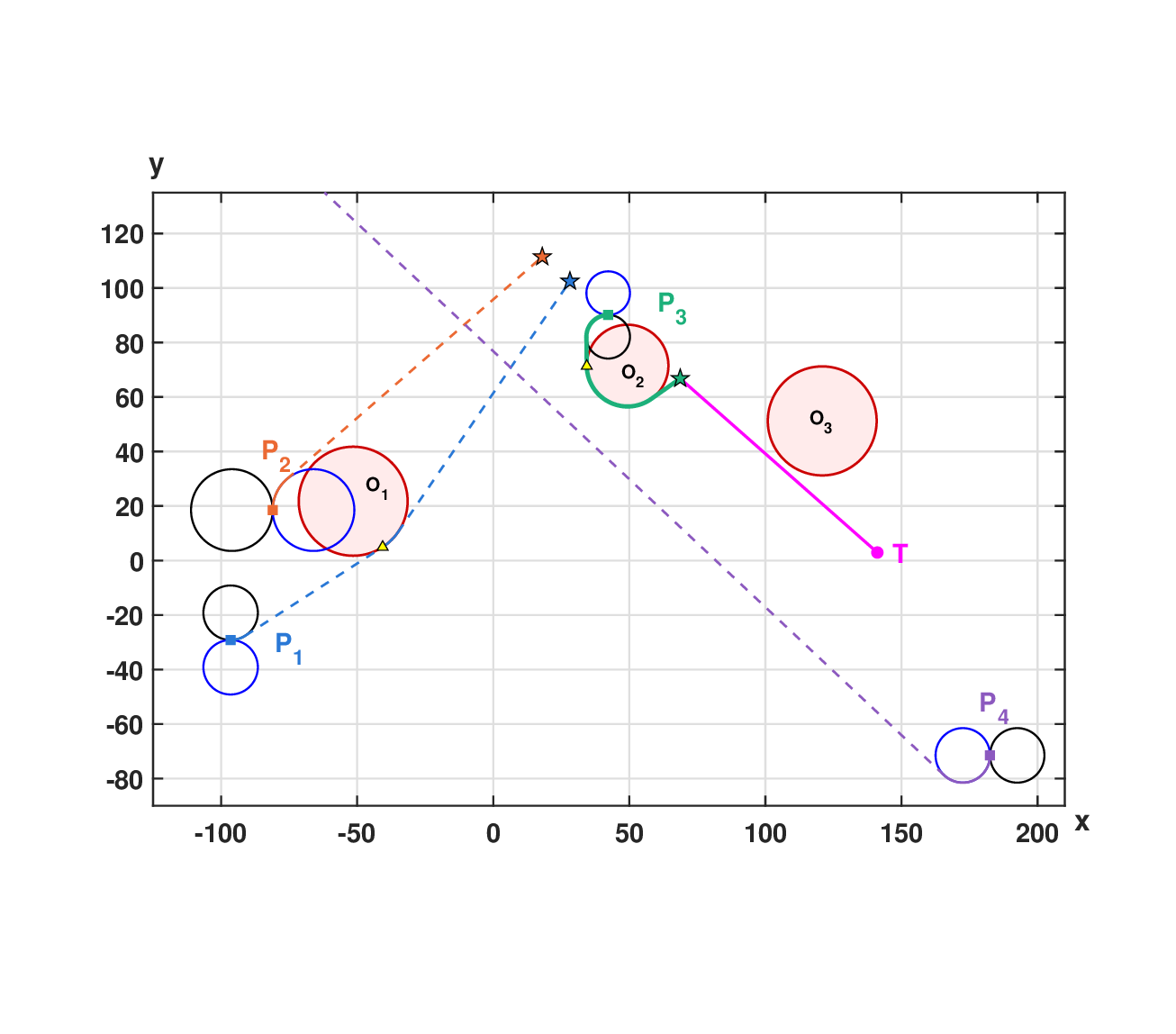}\vspace{-12mm}
\includegraphics[width=0.80\textwidth]{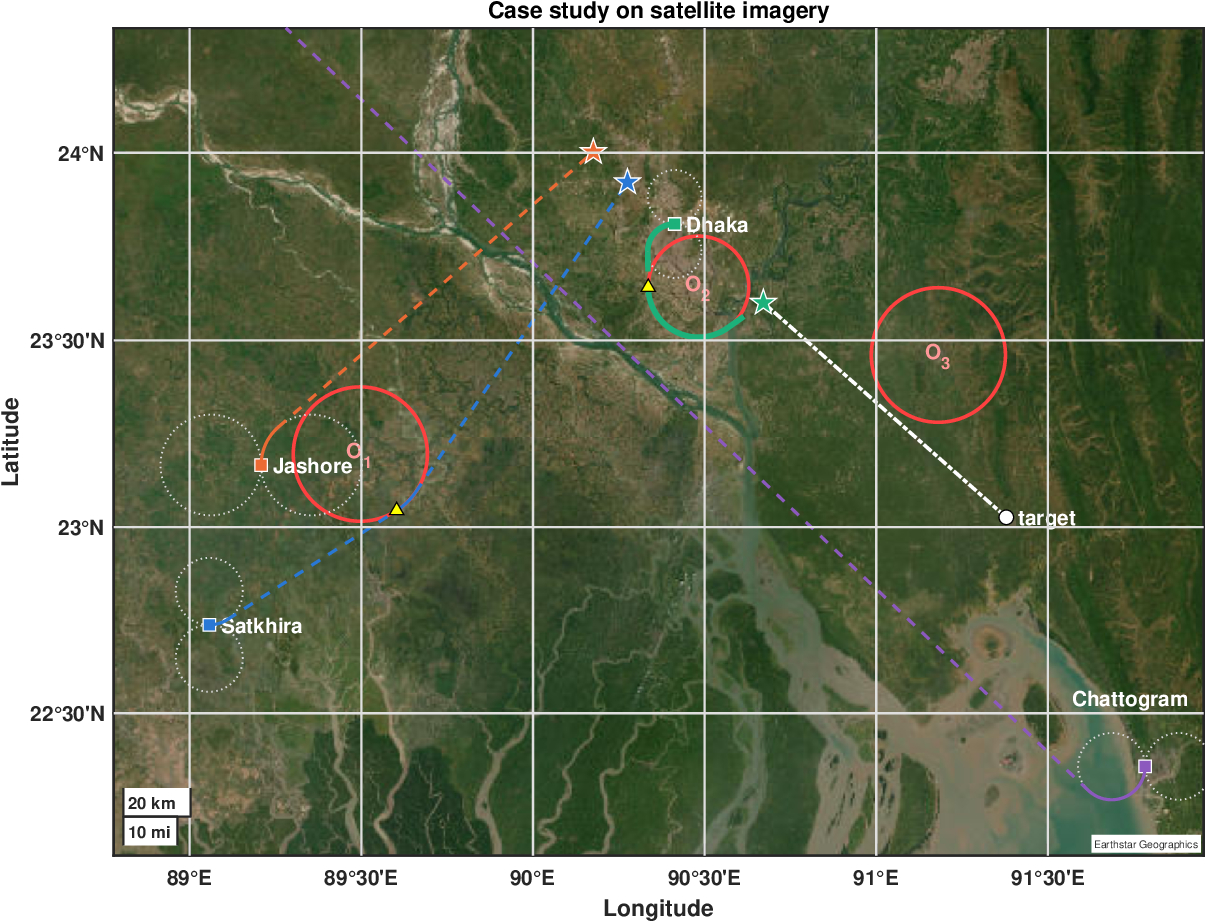}\vspace{5mm}
\caption{\sf\footnotesize Case study (planar coordinates in km): best verified paths of the four vehicles; only the vehicle at Dhaka ($P_3$) is launched. The path of the vehicle at Chattogram ($P_4$) leaves the displayed region; its interception point lies about $1{,}360$ km north-west of Feni.}
\label{fig:case}
\end{figure}

\noindent Table \ref{tab:caseres} and Figure \ref{fig:case} give the results. The vehicle at Dhaka intercepts the intruder after $8.91$ min: it turns left from its westward heading, flies to the boundary of the Narayanganj airspace, follows it counterclockwise for $32.6$ km and meets the intruder at $(68.7,66.8)$, i.e. near $23.60^\circ$N, $90.67^\circ$E, before the intruder reaches the Narayanganj airspace. The airborne vehicle at Satkhira would need $13.90$ min (its direct route crosses the Narail airspace), the fastest vehicle, at Jashore, $15.15$ min, and the nearest vehicle, at Chattogram, more than two hours, since it must chase the intruder from behind. The selected vehicle is thus neither the nearest nor the fastest, and the nearest is by far the worst.

\begin{table}[H]
\caption{\small{\textit{Case study: lower bounds and best verified interception paths (times in minutes). The last column shows the decisions of Algorithm \ref{alg:select}.}}}
\label{tab:caseres}
\footnotesize
\vskip 1em
\centering
\resizebox{\textwidth}{!}{%
\begin{tabular}{|l|c|c|c|c|c|c|}
\hline
Vehicle & Distance [km] & $\underline{T}_k$ & Path type & $L_k$ [km] & $\widehat T_k$ & Algorithm \ref{alg:select} \\
\hline
$P_3$ Dhaka & 131.8 & 7.42 & $LS+LS(O_2)+S_T$ & 69.0 & $\mathbf{8.91}$ & solved (1st), selected \\
$P_1$ Satkhira & 239.8 & 13.61 & $LS+LS(O_1)+S_T$ & 185.3 & 13.90 & pruned \\
$P_2$ Jashore & 222.7 & 15.06 & $RS+S_T$ & 137.2 & 15.15 & pruned \\
$P_4$ Chattogram & 85.2 & 115.84 & $RS+S_T$ & 1454.3 & 125.35 & pruned \\
\hline
\end{tabular}}
\end{table}

\noindent Algorithm \ref{alg:select} solves only the vehicle at Dhaka, which has the smallest lower bound ($7.42$ min). Its verified time, $8.91$ min, is below the lower bounds of the three other vehicles ($13.61$, $15.06$ and $115.84$ min), so all three are pruned and the decision is certified by Proposition \ref{prop:cert}: no other vehicle can intercept the intruder earlier than $8.91$ min, whatever path it flies. The computation takes $81.6$ s instead of $347.4$ s for full enumeration, a reduction of $77\%$, and the bracket on the optimal interception time is $7.42\leq T^*\leq8.91$ min. In this case the pruning is exact even though each vehicle's own problem is solved heuristically: the selected vehicle is optimal among the four, and only the optimality of its path within its own problem (at most $1.49$ min of potential improvement) remains uncertified.
%=======================================================================
\section{Discussion and Limitations}\label{sec:discussion}
%=======================================================================
\noindent \textit{What is certified and what is not.} The lower bound of Lemma \ref{lem:lb}, the pruning rule of Algorithm \ref{alg:select} and the feasibility of every reported path are rigorous, and so is the bracket $T_{\rm LB}\leq T^*\leq T_{\rm feas}$. The value $\widehat T_k$ of one pursuer is not certified to be optimal: it is the best of many verified local solutions over a finite path family. Two pursuers whose computed times differ by less than the uncertainty of the within-pursuer optimisation may therefore be ranked incorrectly; the seed study of Section \ref{sec:seeds} indicates how large this uncertainty is in practice. A certified selection would require a lower bound for each pursuer that accounts for curvature and obstacles, for example from obstacle-free Dubins interception \cite{Meyer2015} or from relaxations of the obstacle-touring problem \cite{Li2026}; tighter bounds would also prune more pursuers.

\noindent \textit{Path family.} The family $\mathcal{F}$ admits $CC$ starts and tours of every ordered subset of obstacles, but it does not contain paths with additional minimum-radius turns between two obstacles or before the final straight segment, nor repeated tours of an obstacle. Under (A3) the common tangent between two obstacle circles is the natural connection, and in the benchmark enlarging the family did not change any pursuer's value (Section \ref{sec:family-num}), but we do not prove that $\mathcal{F}$ contains an optimal path of \eqref{eq:ocp}. The number of path types grows as $N(n,n)$ in \eqref{eq:count}; for many obstacles, the order and subset of toured obstacles should be restricted, for instance to obstacles that intersect the cone between the pursuer and the target region, or chosen by a combinatorial search.

\noindent \textit{Modelling assumptions.} The selection rule relies on (A2) and (A4). If the target can change its heading or speed, each pursuer's interception time becomes uncertain, and launching a second pursuer can reduce the risk that the target escapes; Remark \ref{obs:single} then no longer applies and the selection becomes a robust or stochastic problem. If several targets appear, the interception times $\widehat T_k$ computed here become the costs of an assignment problem \cite{Ahuja2007,Kline2019}. Assumption (A3) excludes obstacles smaller than a pursuer's turning circle; such obstacles can be enlarged to radius $\rho_k$ at the price of conservatism, or require boundary arcs of a different structure \cite{Boissonnat1996}. The planar, constant-speed kinematics do not model wind, altitude changes or vehicle dynamics.

\noindent \textit{Computation.} One pursuer with three obstacles requires about one minute of serial CPU time in our reference implementation. The pursuers are independent, so Step 2 of Algorithm \ref{alg:select} can be parallelised across vehicles, or run on board each vehicle, at the cost of weaker pruning; within one pursuer the path types are also independent. For time-critical use, the number of starting points can be reduced, at the price of a higher risk of missing the best path type of a pursuer.

%=======================================================================
\section{Conclusion}\label{sec:conclusion}
%=======================================================================
\noindent We have presented an obstacle-aware interceptor-selection framework for a team of heterogeneous, curvature-constrained UAVs that respond to a target moving at constant velocity among static circular obstacles, when one vehicle is launched. Each vehicle's interception problem is solved through an arc-parameterised model with launch delays, a path family that includes $CC$ starts and tours of every ordered subset of obstacles, and exact obstacle-avoidance constraints on every segment, and only verified feasible paths are used. A closed-form lower bound and the verified times are combined in a bound-and-select algorithm whose pruning is certified: it returns the same vehicle and time as solving every vehicle, together with a bracket on the optimal team interception time. The experiments showed that the obstacle-aware choice frequently differs from the nearest or the fastest vehicle and from the vehicle with the smallest straight-line time, that the selection changes with the target speed, the launch delay, the turning radius and the obstacle configuration, and that pruning removes a large part of the computation without changing the decision. Future work will address certified per-vehicle bounds, uncertain target motion, several targets and moving obstacles.

%=======================================================================
%\section*{Reproducibility statement}\label{sec:repro}
%=======================================================================
%\noindent All data used in the paper are given in Tables \ref{tab:data} and \ref{tab:casedata} and in the file \texttt{random\_scenarios.csv} (random scenarios of Section \ref{sec:compare}). The accompanying MATLAB package implements exactly the model \eqref{eq:nlp}, the family $\mathcal{F}$, the verification tests, the two-phase multistart protocol and Algorithm \ref{alg:select}, with the settings of Section \ref{sec:setup}; the scripts \texttt{run\_benchmark}, \texttt{run\_seeds}, \texttt{run\_sensitivity}, \texttt{run\_random\_comparison} and \texttt{run\_casestudy} reproduce the tables and figures of Sections \ref{sec:numerics} and \ref{sec:case}, and \texttt{checkAssumptions} verifies (A1), (A3) and (A5) before solving. The reference values in the paper were computed with a Python/SciPy implementation of the same model, which is also provided. Because random-number streams and local solvers differ between the two implementations, individual local solutions may differ. As a cross-check, the MATLAB package was run on the benchmark scenario under GNU Octave~8.4, with Octave's \texttt{sqp} solver substituted for \texttt{fmincon}: it reproduced the best path types, the segment lengths and the interception times of Table \ref{tab:results} to four decimals, as well as the pruning decision of Algorithm \ref{alg:select}.

\section{Declarations}
\noindent \textbf{\large {Funding:}} Not applicable.

\noindent \textbf{\large {Declaration of generative AI and AI-assisted technologies in the writing process:}} During the preparation of this manuscript, the authors used ChatGPT (OpenAI) to assist with language editing, grammatical corrections and improvements in the clarity and presentation of the text. The authors reviewed and edited the content and take full responsibility for the final manuscript. 

\noindent \textbf{\large {Conflicts of interests:}} The authors declare no conflict of interest.

\end{document}